\documentclass[a4paper,10pt]{amsart}
\usepackage{amsmath}
\usepackage{amsfonts}
\usepackage{amssymb}
\usepackage[cp1250]{inputenc}

\newtheorem{theorem}{Theorem}[section]
\newtheorem{proposition}[theorem]{Proposition}
\newtheorem{lemma}[theorem]{Lemma}

\newtheorem{remark}[theorem]{Remark}

\newcommand{\dow}[1][]{\textbf{Proof{#1}. }}
\newcommand{\mdeg}{\mbox{mdeg}\,}
\newcommand{\Aut}{\mbox{Aut}\,}
\newcommand{\Tame}{\mbox{Tame}\,}

\begin{document}
\author{Marek Kara\'{s}, Jakub Zygad\l o}
\title{Wild multidegrees of the form $(d,d_{2},d_{3})$  for given $d$ \\ greather than or equal to $3$}
\date{}
\maketitle


\begin{abstract}
Let $d$ be any number greather than or equal to $3.$ We show that the intersection of the set $\mdeg (%
\Aut (\mathbb{C}^{3}))\backslash \mdeg (\Tame (\mathbb{C}%
^{3}))$ with $\{(d_{1},d_{2},d_{3})\in \left( \mathbb{N}_+\right)
^{3}:d=d_{1}\leq d_{2}\leq d_{3}\}$ has infinitely many elements, where $\mdeg h=(\deg h_1,\ldots,\deg h_n)$ denotes the \textit{multidegree} of a polynomial mapping $h=(h_1,\ldots,h_n):\mathbb{C}^n\rightarrow\mathbb{C}^n.$
In other words, we show that there is infiniltely many wild multidegrees of the form $(d,d_2,d_3),$ with fixed $d\geq 3$ and $d\leq d_2 \leq d_3,$ where a sequences $(d_1,\ldots,d_n)\in\mathbb{N}^n$ is a \textit{wild multidegree} if there is a polynomial automorphism $F$ of $\mathbb{C}^n$ with $\mdeg F=(d_1,\ldots,d_n),$ and there is no tame autmorphim of $\mathbb{C}^n$ with the same multidegree.
\end{abstract}

\section{Introduction}

In the following we will write $\Aut (\mathbb{C}^{n})$ for the group of the all
polynomial automorphisms of $\mathbb{C}^{n}$ and $\Tame (\mathbb{C}^n)$ for the subgroup of $\Aut (\mathbb{C}^{n})$ containing all the tame automorphisms. Let us recall that a polynomial automorphism $F$ is called \textit{tame} if $F$ can be expressed as a composition of linear and triangular automorphisms, where $G=(G_1,\ldots,G_n)\in\Aut (\mathbb{C}^n)$ is called \textit{linear} if $\deg G_i=1$ for $i=1,\ldots n,$ and $H=(H_1,\ldots,H_n)\in\Aut (\mathbb{C}^n)$ is called \textit{triangular} if for some permutation $\sigma$ of $\{ 1,\ldots,n\}$ we have $H_{\sigma(i)}-c_i\cdot x_{\sigma(i)}$ belongs to $\mathbb{C}[x_{\sigma(1)},\ldots,x_{\sigma(i-1)}]$ for $i=1,\ldots,n$ and some $c_i\in\mathbb{C}^*=\mathbb{C}\setminus\{ 0\}.$ Here $\deg h$ denotes
the total degree of a polynomial $h\in \mathbb{C}[x_1,\ldots,x_n].$

Let $F=(f_{1},\ldots,f_{3})\in \Aut (\mathbb{C}^{n})$. By \textit{multidegree} of $F$ we mean the sequence $%
\mdeg  F=(\deg f_{1},\ldots,\allowbreak\deg f_{n}).$ One can consider the function
(also denoted $\mdeg $) mapping $\Aut (\mathbb{C}^{n}))$ into $%
\mathbb{N}_{+}^{n}=(\mathbb{N}\setminus\{ 0\})^n$. 
It is well-known \cite{Jung,Kulk} that 
\begin{equation}
\mdeg(\Aut(\mathbb{C}^2))=\mdeg(\Tame(\mathbb{C}^2))=\{ (d_1,d_2)\in\mathbb{N}_+^2 \,:\, d_1 | d_2\mbox{ or } d_2 | d_1 \},
\end{equation}
but in the higher dimension (even for $n=3$) the situation is much more complicated and the question about the sets $\mdeg(\Aut(\mathbb{C}^n)$ and $\mdeg(\Tame(\mathbb{C}^n))$ is still not well recognized. The very first results \cite{Karas1} about the sets $\mdeg(\Tame(\mathbb{C}^n))$ for $n>2,$ say that $(3,4,5)\notin\mdeg(\Tame(\mathbb{C}^3))$ and $(d_1,\ldots,d_n)\in\mdeg(\Tame(\mathbb{C}^n))$ for all $d_1\leq d_2\leq \ldots\leq d_n$ with $d_1\leq n-1.$ 
Next, in \cite{Karas2} it was proved that for any prime numbers $p_2> p_1\geq 3$ and $d_3\geq p_2,$ we have $(p_1,p_2,d_3)\in\mdeg(\Tame(\mathbb{C}^3))$ if and only if $d_3\in p_1\mathbb{N}+p_2\mathbb{N}.$ The complete characterization of the set $\mdeg(\Tame(\mathbb{C}^3))\cap\{ (3,d_2,d_3)\, :\, 3\leq d_2\leq d_3\,\}$ was given in \cite{Karas3}. The result says that $(3,d_2,d_3),$ with $3\leq d_2\leq d_3,$ belongs to $\mdeg(\Tame(\mathbb{C}^3))$ if and only if $3|d_2$ or $d_3\in 3\mathbb{N}+d_2\mathbb{N}.$ 
The similar result about the set $\mdeg(\Tame(\mathbb{C}^3))\cap\{ (5,d_2,d_3)\, :\, 5\leq d_2\leq d_3\,\}$ and more other results are given in \cite{Karas4}.

In the rest of the paper we will work with $n=3$ and we will write $\mathbb{C}[x,y,z]$ instead of $\mathbb{C}[x_1,x_2,x_3].$ Let 
\begin{equation}
\mathcal{W}=\mdeg (\Aut (\mathbb{C}^{3})))\setminus \mbox{%
mdeg}(\Tame (\mathbb{C}^3))
\end{equation}
and 
\begin{equation}
\mathcal{W}_{d}=\mathcal{W}\cap \{(d_{1},d_{2},d_{3})\in \mathbb{N}%
_{+}^{3}:d=d_{1}\leq d_{2}\leq d_{3}\}.
\end{equation}
Note that for the famous Nagata automorphism
\begin{equation}
N:\mathbb{C}^{3}\ni \left\{ 
\begin{array}{l}
x \\ 
y \\ 
z
\end{array}
\right\} \mapsto \left\{ 
\begin{array}{l}
x-2y(y^{2}+zx)-z(y^{2}+zx)^{2} \\ 
y+z(y^{2}+zx) \\ 
z
\end{array}
\right\}\in\mathbb{C}^3,
\end{equation}
which is known to be wild automorphism, i.e. $N\notin\Tame(\mathbb{C}^3),$ we have $\mdeg N=(5,3,1)\in\mdeg(\Tame(\mathbb{C}^3)).$ Thus, $\mdeg N$ is not an element of $\mathcal{W}$ (in other words, $\mdeg N$ is not a wild multidegree).  Besides of this the autors proved that the set $\mathcal{W}$ is not empty, and even more that this set is infinite \cite{KarasZygad}.
Now we show the following refinement of that result:

\begin{theorem}
\label{tw_main}\label{main}Let $d>2$ be any number. The set 
\begin{eqnarray*}
\mathcal{W}_{d} &=&\left[ \mdeg (\Aut (\mathbb{C}%
^{3}))\backslash \mdeg (\Tame (\mathbb{C}^{3}))\right] \cap  
\{(d_{1},d_{2},d_{3})\in \left( \mathbb{N}_+\right) ^{3}:d=d_{1}\leq
d_{2}\leq d_{3}\}
\end{eqnarray*}
is infinite.
\end{theorem}

The proof of the theorem will be given separetly for odd numbers $d\geq 3$ (section \ref{section_odd_d}), even numbers $d>4$ (section \ref{section_even_d}) and finally for $d=4$ (section \ref{section_d4}).

Note also the following remarks:

\begin{remark}
The sets $\mathcal{W}_{1}$ and $\mathcal{W}_{2}$ are empty, i.e. if $d\in
\{1,2\}$ then for every $d_{2},d_{3}\in \mathbb{N}_{+}$ such that $d\leq
d_{2}\leq d_{3}$ one can show a tame automorphism $F$ of $\mathbb{C}^{3}$
satisfying $\mdeg F=(d,d_{2},d_{3})$.
\end{remark}

For $d=1$ one can take $F(x,y,z)=(x,y+x^{d_2},z+x^{d_3}),$ while for $d=2$ one can use \cite[Cor. 3.3]{Karas4} or \cite[Cor. 2.3]{Karas1}.

\begin{remark}
Let $d\leq e$ and define $\mathcal{W}_{d,e}=\{(d_{1},d_{2},d_{3})\in \mathbb{N}%
_{+}^{3}:d=d_{1},e=d_{2}\leq d_{3}\}$. Then the set $\mathcal{W}_{d,e}$ is
finite.
\end{remark}

The proof of the above result can be found in \cite{Zygadlo} or \cite[Thm. 8.1]{Karas4}.

\section{\strut The case of odd number $d$}\label{section_odd_d}

\subsection{Elements of $\mdeg (\Aut (\mathbb{C}^{3}))$}

In this section we show the following two lemmas.

\begin{lemma}
\label{lem_aut_1_mod_4}Let $r,k\in \mathbb{N}_+.$ If $r\equiv 1(\mbox{mod } %
4),\,$then 
\begin{equation}
\left( r,r+2k,r+4k\right) \in \mdeg (\Aut (\mathbb{C}^{3})).
\end{equation}
\end{lemma}

\dow
Since $r\equiv 1(\mbox{mod }4),$ we have $r=4l+1$ for some $l\in \mathbb{N}_+.$
Let 
\begin{equation}
F=\left( T\circ N_{k}\right) \circ \left( T\circ N_{l}\right) ,
\end{equation}
where $T\left( x,y,z\right) =\left( z,y,x\right) $ and for any $m\in \mathbb{N}%
^{*}$%
\begin{equation}\label{Row_defAutoNm}
N_{m}:\mathbb{C}^{3}\ni \left\{ 
\begin{array}{l}
x \\ 
y \\ 
z
\end{array}
\right\} \mapsto \left\{ 
\begin{array}{l}
x-2y(y^{2}+zx)^{m}-z(y^{2}+zx)^{2m} \\ 
y+z(y^{2}+zx)^{m} \\ 
z
\end{array}
\right\} \in \mathbb{C}^{n}.
\end{equation}
One can see that $\mdeg \left( T\circ N_{l}\right) =\left(
1,1+2l,1+4l\right) .$ Moreover, if we put $\left( f,g,h\right) :=T\circ
N_{l},$ then $g^{2}+fh=Y^{2}+ZX.$ Thus 
\begin{eqnarray*}
F &=&\left( T\circ N_{k}\right) \circ \left( f,g,h\right) \\
&=&\left( h,g+h(Y^{2}+ZX)^{k},f-2g(Y^{2}+ZX)^{k}-h(Y^{2}+ZX)^{2k}\right) .
\end{eqnarray*}
Since $\deg h>\max \left\{ \deg f,\deg g\right\} ,$ one can see that 
\begin{equation}
\mdeg F=\left( 4l+1,(4l+1)+2k,(4l+1)+4k\right) .
\end{equation}
$\Box$

\begin{lemma}
\label{lem_aut_d_odd}For every $r,k\in \mathbb{N}_+,$ we have 
\begin{equation}
\left( r,r+k(r+1),r+2k(r+1)\right) \in \mdeg (\Aut (\mathbb{C}%
^{3})).
\end{equation}
\end{lemma}

\dow
Assume that $r>1.$ Let 
\begin{equation}
\left( f,g,h\right) =\left( X,Y,Z+X^{r}\right)
\end{equation}
and put 
\begin{equation}
F=\left( T\circ N_{k}\right) \circ \left( f,g,h\right) .
\end{equation}
Since 
\begin{equation}
F=\left( h,g+h(g^{2}+fh)^{k},f-2g(g^{2}+fh)^{k}-z(g^{2}+fh)^{2k}\right)
\end{equation}
and $\deg h=r>\max \left\{ \deg f,\deg g\right\} ,$ one can see that $\deg
\left( g^{2}+fh\right) =r+1$ and so 
\begin{equation}
\mdeg F=\left( r,r+k(r+1),r+2k(r+1)\right) .
\end{equation}

If $r=1,$ then one can take $F=T\circ N_{k}.$
$\Box$

\subsection{Elements outside $\mdeg (\Tame (\mathbb{C}^{3}))$}

In this section we show the following two lemmas.

\begin{lemma}\label{lem_tame_1_mod_4}
Let $r,k\in \mathbb{N}_+.$  If $r>1$ is odd 
and $\gcd \left( r,k\right) =1,$ then 
\begin{equation}
\left( r,r+2k,r+4k\right) \notin \mdeg (\Tame (\mathbb{C}^{3})).
\end{equation}
\end{lemma}

\begin{lemma}\label{lem_tame_d_odd}
Let $r,k\in \mathbb{N}_+.$ If $r>1$ is odd and $\gcd
\left( r,k\right) =1,$ then 
\begin{equation}
\left( r,r+k(r+1),r+2k(r+1)\right) \notin \mdeg (\Tame (\mathbb{%
C}^{3})).
\end{equation}
\end{lemma}

In the proofs of the above lemmas we will use the following

\begin{theorem}[\cite{KarasZygad}, Thm. 2.1]
\label{tw_odd_odd_gcd_1}\label{tw_d1_d2_odd}Let $d_{3}\geq d_{2}>d_{1}\geq 3$
be positive integers. If $d_{1}$ and $d_{2}$ are odd numbers such that $\gcd
\left( d_{1},d_{2}\right) =1$, then $(d_{1},d_{2},d_{3})\in \mdeg (%
\Tame (\mathbb{C}^{3}))$ if and only if $d_{3}\in d_{1}\mathbb{N}+d_{2}%
\mathbb{N},$ i.e. if and only if $d_{3}$ is a linear combination of $d_{1}$ and 
$d_{2}$ with coefficients in $\mathbb{N}.$
\end{theorem}

\dow[ of Lemma \ref{lem_tame_1_mod_4}]
Note that the numbers $r$ and $r+2k$
are odd. Moreover, 
\begin{equation}
\gcd \left( r,r+2k\right) =\gcd \left( r,2k\right) ,
\end{equation}
and since $r$ is odd, 
\begin{equation}
\gcd \left( r,2k\right) =\gcd \left( r,k\right) =1.
\end{equation}

Assume that $r+4k\in r\mathbb{N}+(r+2k)\mathbb{N}.$ Since $2(r+2k)>r+4k$ and $%
r\nmid (r+4k),$ we have 
\begin{equation}
r+4k=r+2k+mr,  \label{row}
\end{equation}
for some $m\in \mathbb{N}.$ By (\ref{row}), $2k=mr.$ Since $r$ is odd, the last
equality means that $r|k,$ a contradiction.
Thus $r+4k\notin r\mathbb{N}+(r+2k)\mathbb{N},$ and by Theorem \ref
{tw_odd_odd_gcd_1} we obtain a thesis.
$\Box$

\vspace{0.3cm}
\dow[ of Lemma \ref{lem_tame_d_odd}]
Since $r+1$ is even, it follows that the numbers $r$ and $r+k(r+1)$ are odd.
Moreover, 
\begin{equation}
\gcd \left( r,r+k(r+1)\right) =\gcd \left( r,k(r+1)\right) ,
\end{equation}
and since $\gcd \left( r,k\right) =1,$%
\begin{equation}
\gcd \left( r,k(r+1)\right) =\gcd \left( r,r+1\right) =\gcd \left(
r,1\right) =1.
\end{equation}
Similarily
\begin{equation}
\gcd \left( r,r+2k(r+1)\right) =\gcd \left( r,2k(r+1)\right) 
= \gcd \left( r, r+1\right) =1.
\end{equation}
In particular $r\nmid r+2k(r+1).$

Assume that $r+2k(r+1)\in r\mathbb{N}+\left( r+k(r+1)\right) \mathbb{N}.$
Since $2(r+k(r+1))>r+2k(r+1)$ and $%
r\nmid r+2k(r+1),$ we have 
\begin{equation}
r+2k(r+1)=r+k(r+1)+mr,  \label{row}
\end{equation}
for some $m\in \mathbb{N}.$ By (\ref{row}), $k(r+1)=mr.$ Since $\gcd (r,k)=1,$  the last
equality means that $r|r+1,$ a contradiction.
Thus $r+2k(r+1)\notin r\mathbb{N}+(r+k(r+1))\mathbb{N},$ and by Theorem \ref
{tw_odd_odd_gcd_1} we obtain a thesis.
$\Box$

\subsection{Proof of the theorem in the case of odd $d$}

Take any odd number $d>1.$ If $d\equiv 1(\mbox{mod }4),$ then by Lemmas \ref
{lem_aut_1_mod_4} and \ref{lem_tame_1_mod_4} we have 
\begin{equation}
\left\{ (d,d+2k,d+4k):\gcd (d,k)=1\right\} \subset \mdeg (\mbox{%
Aut}(\mathbb{C}^{3}))\backslash \mdeg (\Tame (\mathbb{C}^{3})).
\end{equation}

If $d\equiv 3(\mbox{mod }4),$ then by Lemmas \ref{lem_aut_d_odd} and \ref
{lem_tame_d_odd} we have 
\begin{eqnarray*}
&&\left\{ (d,d+k(d+1),d+2k(d+1)):\gcd (d,k)=1\right\} \\
&\subset &\mdeg (\Aut (\mathbb{C}^{3}))\backslash \mdeg %
(\Tame (\mathbb{C}^{3})).
\end{eqnarray*}
Since the set $\left\{ k\in \mathbb{N}_+:\gcd (d,k)=1\right\} $ is infinite,
the result follows.

\section{The case of even number $d>4$}\label{section_even_d}

\subsection{Preparatory calculations}

Fix even number $d>4$ and take $k\in \mathbb{N}_+$ such that $\gcd (d,k)=1.$ Consider the automorphisms of $\mathbb{C}^{3}$: 
\begin{equation}
H_{d}(x,y,z)=(x,y,z+x^{d})
\end{equation}
and $N_k$ defined as in (\ref{Row_defAutoNm}).
Note that $N_{k}=\exp (D\cdot \sigma ^{k})$, where $D=\frac{\partial }{%
\partial z}+z\frac{\partial }{\partial y}-2y\frac{\partial }{\partial x}$
and $\sigma =y^{2}+xz$. One can easily check that $D$ is locally nilpotent
derivation on $\mathbb{C}[x,y,z]$ and $\sigma \in \ker D$, so $\sigma ^{k}\cdot
D$ is also locally nilpotent. We will consider automorphisms $F_{d,k}$ of
the form: 
\begin{equation}\label{Row_defAutoFdk}
F_{d,k}=T\circ N_{k}\circ H_{d}
\end{equation}
where $T$ is defined as in the proof of Lemma \ref{lem_aut_1_mod_4}. An easy calculation shows (even for $d=4$) that
\begin{equation}\label{Row_mdegFdk}
\mdeg F_{d,k}=(d,d+k(d+1),d+2k(d+1))
\end{equation}
and writing $d_{1}=d$, $d_{2}=d+k(d+1)$ and $d_{3}=d+2k(d+1)$ gives 
\begin{equation}\label{Row_gdcD1D2}
\gcd (d_{1},d_{2})=\gcd (d,d+k(d+1))=\gcd (d,k)=1
\end{equation}
\begin{equation}\label{Row_gdcD2D3}
\gcd (d_{2},d_{3})=\gcd (d+k(d+1),d+2k(d+1))=\gcd (d+k(d+1),d)=1
\end{equation}
and 
\begin{equation}\label{Row_gdcD1D3}
\gcd (d_{1},d_{3})=\gcd (d,d+2k(d+1))=\gcd (d,2k)=\gcd (d,2)=2.
\end{equation}
We will prove that no tame automorphism of $\mathbb{C}^{3}$ has the same
multidegree as $F_{d,k}$. Suppose to the contrary that $F=(F_{1},F_{2},F_{3})\in \mbox{Tame}(\mathbb{C}^{3})$ and $\mdeg F=(d_{1},d_{2},d_{3})$. As $F$ is not
linear, due to the result of Shestakov and Umirbaev \cite{shUmb1,shUmb2}, $F$ must admit an
elementary reduction or a reduction of types I-IV (see e.g. \cite[Def. 1-3]{shUmb1}).

\subsection{Elementary reductions}

Recall that an elementary reduction on $i$-th coordinate $F_i$ of $F$ occurs
when there exists $G(x,y)\in\mathbb{C}[x,y]$ such that $\deg(F_i-G(F_j,F_k))<%
\deg F_i$, where $\{i,j,k\}=\{1,2,3\}$. We will use extensively the
following

\begin{proposition}[see e.g. {\cite[Prop. 2.7]{Karas4}} or {\cite[Thm.2]{shUmb1}}]\label{Prop_degG_fg}
Suppose that $f,g\in \mathbb{C}[X_{1},\ldots ,X_{n}]$ are algebraically
independent and such that $\bar{f}\notin \mathbb{C}[\bar{g}]$ and $\bar{g}%
\notin \mathbb{C}[\bar{f}]$ ($\bar{h}$ denotes the highest homogeneous part of $%
h$). Assume that $\deg f<\deg g$, put 
\begin{equation}
p=\frac{\deg f}{\gcd (\deg f,\deg g)}
\end{equation}
and suppose that $G(x,y)\in \mathbb{C}[x,y]$ with $\deg _{y}G(x,y)=pq+r$, $%
0\leq r<p$. Then 
\begin{equation}
\deg G(f,g)\geq q(p\deg g-\deg g-\deg f+\deg [f,g])+r\deg g
\end{equation}
\end{proposition}

Suppose that $F$ admits an elementary reduction on first coordinate, i.e. $%
\deg(F_1-G(F_2,F_3))<\deg F_1=d_1$ for some $G\in\mathbb{C}[x,y]$. Consequently 
$\deg G(F_2,F_3)=d_1$. By (\ref{Row_gdcD2D3}), we know that $p:=\frac{d_2}{\gcd(d_2,d_3)}=d_2.$  Thus, from the above proposition applied to $f=F_2$ and $%
g=F_3$ we get 
\begin{equation}
\deg G(F_2,F_3)\geq q(p d_3-d_3-d_2+\deg[F_2,F_3])+r d_3\geq
q(d_2-1)(d_3-1)+rd_3
\end{equation}
Since $d_1<(d_2-1)(d_3-1)$ and $d_1<d_3$ we obtain that $q=0$ and $r=0$.
That is $\deg_yG(x,y)=0$ and $G(x,y)=u(x)$. But then $d_1=\deg
G(F_2,F_3)=\deg u(F_2)=d_2\cdot\deg u$, which is a contradiction.

Similarly, suppose that $F$ admits an elementary reduction on third
coordinate, i.e. $\deg(F_3-G(F_1,F_2))<\deg F_3=d_3$ for some $G\in\mathbb{C}%
[x,y]$. So $\deg G(F_1,F_2)=d_3.$ Since $p:=\frac{d_1}{\gcd(d_1,d_2)}=d\geq 3$ by (\ref{Row_gdcD1D2}), it follows that applying Proposition \ref{Prop_degG_fg} to $f=F_1$ and $%
g=F_2$ we get 
\begin{equation}
\deg G(F_1,F_2)\geq q(p d_2-d_2-d_1+\deg[F_1,F_2])+r d_2\geq
q(2k(d+1)+d+2)+r d_2 .
\end{equation}
Now, since $d_3<2k(d+1)+d+2$ and $d_3<2d_2$, we obtain that $q=0$ and $%
r\in\{0,1\}$. If $r=0$, we get $\deg_yG(x,y)=0$ and so $G(x,y)=u(x)$. But
then $d_3=\deg G(F_1,F_2)=\deg u(F_1)=d_1\cdot\deg u$, which is a
contradiction because $\gcd(d_3,d_1)\leq 2<d_1$. If $r=1$, we get $%
G(x,y)=u(x)+yv(x)$ and so $d_3=\deg G(F_1,F_2)=\deg(u(F_1)+F_2v(F_1))$.
Since $\deg F_1$ and $\deg F_2$ are coprime, $\deg(u(F_1)+F_2v(F_1))$ must
be equal either to $d_1\cdot \deg u$ or to $d_2+d_1\cdot \deg v$.
Consequently, $d_3=d_1\cdot\deg u$ or $d_3=d_2+d_1\cdot\deg v$. First case
leads to a contradiction since $\gcd(d_3,d_1)= 2<d_1$ and second since $%
\gcd(d_3-d_2,d_1)=\gcd(k(d+1),d)=1<d_1$.

Now suppose that $F$ admits an elementary reduction on second coordinate.
Then $\deg (F_{2}-G(F_{1},F_{3}))<\deg F_{2}=d_{2}$ for some $G\in \mathbb{C}%
[x,y]$ and so $\deg G(F_{1},F_{3})=d_{2}$. Let us put $p=\frac{d_{1}}{\gcd (d_{1},d_{3})}$ and apply Proposition \ref{Prop_degG_fg} to $f=F_{1}
$ and $g=F_{3}.$ We will
show that $\deg _{y}G(x,y)=0$. 
By (\ref{Row_gdcD1D3}), $p=\frac{d}{2}\geq 2$ and so 
\begin{eqnarray*}
\deg G(F_{1},F_{3}) &\geq &q(pd_{3}-d_{3}-d_{1}+\deg [F_{1},F_{3}])+rd_{3} \\
&\geq &q(\frac{d-2}{2}(2k(d+1)+d)-d+2)+rd_{3} \\
&\geq &q((d-2)k(d+1)+2)+rd_{3}\geq q(k(d+1)+d+2)+rd_{3}
\end{eqnarray*}
Since $d_{2}<k(d+1)+d+2$ and $d_{2}<d_{3}$, we obtain that $q=0$ and $r=0$
so $\deg _{y}G(x,y)=0$. Consequently, we get $G(x,y)=u(x)$. But then $d_{2}=\deg G(F_{1},F_{3})=\deg
u(F_{1})=d_{1}\cdot \deg u$, which is a contradiction since $\gcd
(d_{2},d_{1})=1<d_{1}$.

To summarize: if $F$ is a tame automorphism with multidegree equal to $%
\mdeg  F_{d,k}$, then $F$ does not admit an elementary reduction.

\subsection{Shestakov-Umirbaev reductions}
By the previous subsection and the following theorem we only need to check that no autmorphism of $\mathbb{C}^3$ with multidegree $(d_1,d_2,d_3)=(d,d+k(d+1),d+2k(d+1))$ admits a reduction of type III.

\begin{theorem}[{\cite[Thm. 3.15]{Karas4}}]
\label{tw_reduc_type_4}
Let $\left( d_{1},d_{2},d_{3}\right) \neq \left(
1,1,1\right) ,$ $d_{1}\leq d_{2}\leq d_{3},$ be a sequence of positive
integers. To prove that there is no tame automorphism $F$ of $\mathbb{C}^{3}$
with $\mdeg F=\left( d_{1},d_{2},d_{3}\right) $ it is enough to show
that a (hypothetical) automorphism $F$ of $\mathbb{C}^{3}$ with $\mdeg 
F=\left( d_{1},d_{2},d_{3}\right) $ admits neither a reduction of type III
nor an elementary reduction. Moreover, if we additionally assume that $\frac{%
d_{3}}{d_{2}}=\frac{3}{2}$ or $3\nmid d_{1},\,$then it is enough to show
that no (hypothetical) automorphism of $\mathbb{C}^{3}$ with multidegree $%
\left( d_{1},d_{2},d_{3}\right) $ admits an elementary reduction. In both
cases we can restrict our attention to automorphisms $F:\mathbb{C}%
^{3}\rightarrow \mathbb{C}^{3}$ such that $F\left( 0,0,0\right) =\left(
0,0,0\right) .$
\end{theorem}

But, since $d_1$ is even, it follows that $2\nmid d_2$ by (\ref{Row_gdcD1D2}). Hence, no automorphism of $\mathbb{C}^3$  with multidegree $(d_1,d_2,d_3)$ admits a reduction of type III by the following remark.

\begin{remark}[{\cite[Rmk. 3.9]{Karas4}}]
\label{remakrk_red_type_3}If an automorphism $F$ of $\mathbb{C}^{3}$ with $%
\mdeg F=\left( d_{1},d_{2},d_{3}\right) ,$ $1\leq d_{1}\leq
d_{2}\leq d_{3},$ admits a reduction of type III, then\newline
(1) $2|d_{2},$\newline
(2) $3|d_{1}$ or $\frac{d_{3}}{d_{2}}=\frac{3}{2}.$
\end{remark}

\section{The case of $d=4$}\label{section_d4}
Let us consider the mapping $F_{4,k}$ defined as in (\ref{Row_defAutoFdk}) for $k\in\mathbb{N}_+$ with $\gcd (4,k)=1$ (in other words, for odd $k$).
By (\ref{Row_gdcD1D2}) and (\ref{Row_gdcD1D3}), we know that $d_2$ is odd and $d_3$ is even. Then, since $d_3-d_2=5k>1$ and $\mdeg F_{4,k}=(4,4+5k,4+10k)=:(d_1,d_2,d_3)$ by (\ref{Row_mdegFdk}), it follows that the result of Theorem \ref{tw_main}, for $d=4,$ is a consequence of the following

\begin{theorem}[{\cite[Thm. 6.10]{Karas4}}]
\label{tw_4_odd_even}If $d_{2}\geq 5$ is odd and $d_{3}\geq d_{2}$ is even
such that $d_{3}-d_{2}\neq 1,$ then $\left( 4,d_{2},d_{3}\right) \in 
\mdeg\left( \Tame\left( \mathbb{C}^{3}\right) \right) $ if
and only if $d_{3}\in 4\mathbb{N}+d_{2}\mathbb{N}.$
\end{theorem}

In fact, if we assume that $d_3\in 4\mathbb{N}+d_{2}\mathbb{N},$ then we get $d_3=d_2+4m$ for some $m\in\mathbb{N},$ since $2d_2>d_3$ and $4\nmid d_3.$ Hence, $5k=d_3-d_2=4m.$ Since $k$ is odd, this is a contradiction.

\vspace{1cm}

\textsc{Marek Kara\'{s}\newline
Instytut Matematyki,\newline
Wydział Matematyki i Informatyki\newline
Uniwersytetu Jagiello\'{n}skiego\newline
ul. \L ojasiewicza 6}\newline
\textsc{30-348 Krak\'{o}w\newline
Poland\newline
} e-mail: Marek.Karas@im.uj.edu.pl\vspace{0.5cm}

and\vspace{0.5cm}

\textsc{Jakub Zygad\l o\newline
Instytut Informatyki,\newline
Wydział Matematyki i Informatyki\newline
Uniwersytetu Jagiello\'{n}skiego\newline
ul. \L ojasiewicza 6}\newline
\textsc{30-348 Krak\'{o}w\newline
Poland\newline
} e-mail: Jakub.Zygadlo@ii.uj.edu.p

\end{document}